\documentclass[12pt,english]{article}
\usepackage{fullpage}
\usepackage{babel}
\usepackage{exscale,xspace}
\usepackage{amsfonts,amssymb}

\usepackage[
  pdfauthor={Bernd C. Kellner},
  pdfkeywords={Bernoulli numbers, irregular pairs, higher order, index of irregularity,
    strong friendly, Chinese remainder theorem, integer sequence, A092291},
  pdftitle={A conjecture about numerators of Bernoulli numbers related to integer sequence A092291},
  colorlinks,bookmarksnumbered]{hyperref}

\newcommand{\pdiv}{\mid}

\newcommand{\notdiv}{\nmid}
\newcommand{\sep}{:\,\,}

\newcommand{\mod}[1]{({\rm mod\ } #1)}
\def\ord{\mathop{\rm ord}\nolimits}
\def\num{\mathop{\rm num}\nolimits}
\def\lcm{\mathop{\rm lcm}\nolimits}

\newcommand{\QQ}{\mathbb{Q}}

\newcommand{\NN}{\mathbb{N}}

\newcommand{\refeqn}[1]{(\ref{#1})}

\newcommand{\qed}{\parbox{0cm}{}\hspace*{\fill} $\Box$}
\newcommand{\eulerphi}{\varphi}

\newcommand{\IRR}{\Psi^{\rm irr}}
\newcommand{\IRP}{\widehat{\Psi}^{\rm irr}}

\newtheorem{prop}{Proposition}[section]
\newtheorem{theorem}[prop]{Theorem}
\newtheorem{corl}[prop]{Corollary}
\newtheorem{lemma}[prop]{Lemma}
\newtheorem{conj}[prop]{Conjecture}

\newtheorem{thdefin}[prop]{Definition}
\newtheorem{thremark}[prop]{Remark}
\newtheorem{thalgo}[prop]{Algorithm}

\newenvironment{defin}{\begin{thdefin}\rm}{\end{thdefin}}
\newenvironment{remark}{\begin{thremark}\rm}{\end{thremark}}
\newenvironment{algo}{\begin{thalgo}\rm}{\end{thalgo}}

\newenvironment{proof}{\medskip\noindent%
\textsc{Proof.}}{\qed\medskip}

\newenvironment{proofof}[1]{\medskip\noindent%
\textsc{Proof of #1.}}{\qed\medskip}

\begin{document}

\nocite{*}

\title{A conjecture about numerators of Bernoulli numbers related to Integer Sequence A092291}
\author{Bernd C. Kellner}
\date{}
\maketitle
\vspace*{-0.8cm}

\abstract{In this paper we disprove a conjecture about numerators of divided  Bernoulli numbers
$B_n/n$ and $B_n/n(n-1)$ which was suggested by Roland Bacher. We give some counterexamples.
Finally, we extend the results to the general case.}
\smallskip

\textbf{Keywords:} Bernoulli number, Kummer congruences, irregular pair,
Chinese remainder theorem
\smallskip

\textbf{Mathematics Subject Classification 2000:} 11B68

\parindent 0cm

\section{Introduction}

Let $B_n$ be the $n$-th Bernoulli number with $n \geq 0$. They are defined by
the power series
\[
  \frac{z}{e^z-1} = \sum_{n=0}^\infty B_n \frac{z^n}{n!} \,,
  \qquad |z| < 2 \pi \,,
\]
where all numbers $B_n$ are zero with odd index $n > 1$.
Therefore, we will consider only even indices concerning Bernoulli numbers.
These numbers play an important role in several topics in mathematics. Here,
we are interested in the numbers
\[
    \frac{B_n}{n} \quad \mbox{and} \quad \frac{B_n}{n(n-1)}
\]
which occur, e.g., in approximation formulas of harmonic numbers $H_n$ resp.\
Stirling's approximation of $\log \Gamma(x)$, see \cite[pp.~480--482]{graham94concrete}.
\medskip

Now, we need some basic facts about Bernoulli numbers
which can be found in \cite[Chapter 15]{ireland90classical}.
In 1850 Kummer introduced the following definition.
\begin{defin} \label{def-irrpair}
Let $p$ be an odd prime. A pair $(p,l)$ is called an \textsl{irregular pair} if $p \pdiv B_l$
with $2 \leq l \leq p-3$ and even $l$.
The \textsl{index of irregularity} of $p$ is defined by
\[
    i(p) := \# \{ (p,l) \mbox{\ is an irregular pair} \sep l = 2, 4, \ldots, p-3 \} \,.
\]
Then $p$ is called an \textsl{irregular} prime if $i(p) > 0$, otherwise $p$ is a \textsl{regular}
prime.
\end{defin}
\medskip

Let $\eulerphi$ be the Euler $\eulerphi$-function, then the classical Kummer congruences state
for $n, n'$ even, $p$ prime, and $p-1 \notdiv n$
\begin{equation} \label{eqn-kummer-congr}
   \frac{B_n}{n} \equiv \frac{B_{n'}}{n'} \pmod{p}
\end{equation}
with $n \equiv n' \ \mod{\eulerphi(p)}$.
An easy consequence of the Kummer congruences supplies that the numerator of
$B_n/n$ consists only of irregular primes and that infinitely many irregular primes exist.
Let $(p,l)$ be an irregular pair. Using congruence \refeqn{eqn-kummer-congr} provides
for all $k \in \NN_0$
\begin{equation} \label{eqn-succ-pdiv}
   p \pdiv B_{l+k\eulerphi(p)}/(l+k\eulerphi(p)) \,.
\end{equation}

The following conjecture about numerators of $B_n/n$ and $B_n/n(n-1)$
was suggested by Roland Bacher, see
The On-Line Encyclopedia of Integer Sequences \cite{IntSeq},
Sequence \textbf{A092291}.
First values are given by 574, 1269, 1910, 3384, 1185, 1376, 9611.
The statements will differ by a factor 2,
because we will use only even indices $n$ instead of $2n$.
Define $\num(r)$ as the numerator of a rational number $r$.

\begin{conj} \label{conj-bern-num}
Let $(p,l)$ be an irregular pair with smallest $l$ in case of index of irregularity $i(p) > 1$.
Define
\[
   A(p) = \min_m \left\{ m \,\,\Big|\,\, \num\left( \frac{B_m}{m} \right) \Big/ \
            \num\left( \frac{B_m}{m(m-1)} \right) \ = \ p \right\} .
\]
Then $A(p) = (l-1)p+1$.
\end{conj}

Actually, let $p_n$ be the $n$-th irregular prime, then $A(p_n)/2$ gives
Integer Sequence \textbf{A092291}.

\section{Counterexamples}
\setcounter{equation}{0}

Because the conjecture does not cover all irregular pairs, we will extend our research
to all of them. Note that for example $(157,62)$ and
$(157,110)$ are irregular pairs and the index of irregularity is $i(157)=2$.

\begin{theorem} \label{theor-Ap-small}
Let $(p,l)$ be an irregular pair.
Define
\[
   A(p) = \min_k
     \left\{ m=l+k\eulerphi(p) \,\,\Big|\,\,
        \num\left( \frac{B_m}{m} \right) \Big/ \
           \num\left( \frac{B_m}{m(m-1)} \right) \ = \ p \right\} .
\]
Then $A(p) = (l-1)p+1$ is valid and has smallest possible value
if and only if one of the following cases holds
\begin{enumerate}
\item $l-1$ has no irregular prime factors.
\item If $q$ is an irregular prime divisor of $l-1$, then $q \notdiv B_{(l-1)p+1}/((l-1)p+1)$.
\end{enumerate}
\end{theorem}

\begin{proof}
First of all, we will prove that $A(p)=(l-1)p+1$ is the smallest possible value.
To solve
\[
    \num \left( \frac{B_m}{m} \right) \Big/ \
          \num \left( \frac{B_m}{m(m-1)} \right) \ = \ p \,,
\]
factor $m-1$ must have the form $m-1 = pc$ with some integer $c$ to reduce
the $p$-power of the second numerator. In other words, we must have
\[
    \ord_p \, \num \left( \frac{B_m}{m} \right) = s \quad \mbox{and} \quad
      \ord_p \, \num \left( \frac{B_m}{m(m-1)} \right) = s-1
\]
with some integer $s \geq 1$.
Let $m'$ be the smallest possible value we are searching for.
By \refeqn{eqn-succ-pdiv} we then have
$m' = l + k(p-1)$ and $m'-1 = pc$. This yields
\[
   l-1 + k(p-1) = pc  \quad \mbox{resp.} \quad k \equiv l-1 \pmod{p} \,.
\]
By definition we have $1 < l < p-2$. Thus, $k=l-1$ is the smallest possible value and
finally
\[
   m' = l+(l-1)(p-1) = (l-1)p+1 = A(p) \,.
\]
Now, we have to take care that $m'-1=(l-1)p$ does not delete other irregular prime factors
of the numerator of $B_{m'}/m'$. In case (1) nothing happens. In case (2) an
irregular prime divisor $q$ of $l-1$ must not appear in the numerator of $B_{m'}/m'$.
\end{proof}

Using Kummer congruences \refeqn{eqn-kummer-congr} and property \refeqn{eqn-succ-pdiv}
again, we can now reformulate Conjecture \ref{conj-bern-num} to an extended equivalent
conjecture described only by irregular pairs.

\begin{conj}
Let $(p,l)$ be an irregular pair. If $q$ is an irregular prime divisor of $l-1$
then for all irregular pairs $(q,l')$ the following holds
\[
   (l-1)p \ \not\equiv \ l'-1 \pmod{q-1} \,.
\]
\end{conj}

But this conjecture is \textbf{not} valid.
We have done some calculations for all irregular pairs $(p,l)$ with $p < 1\,000\,000$
using a database of irregular pairs calculated in \cite{irrprime12M}. There are
39\,181 irregular pairs all together, 16\,540 of them have irregular prime divisors of
the corresponding $l-1$ and 149 exceptions occur.
\medskip

The first five exceptions and the last calculated exception
are listed below.
\smallskip

\begin{center}
\begin{tabular}{|c|r|c|c|} \hline
$(p,l)$ & $m= (l-1)p+1$ & $l-1$ & $(q,l')$ \\ \hline\hline
$(6449,4884)$  & $31\,490\,468$ & $19 \cdot 257$ & $(257,164)$ \\ \hline
$(8677,2658)$  & $23\,054\,790$ & $2657$         & $(2657,710)$ \\ \hline
$(11351,1044)$ & $11\,839\,094$ & $7 \cdot 149$  & $(149,130)$ \\ \hline
$(12527,2122)$ & $26\,569\,768$ & $3 \cdot 7 \cdot 101$ & $(101,68)$ \\ \hline
$(15823,482)$  & $7\,610\,864$  & $13 \cdot 37$ & $(37,32)$  \\ \hline
$\ldots$ & \multicolumn{1}{c|}{$\ldots$} & $\ldots$ & \multicolumn{1}{c|}{$\ldots$} \\ \hline
$(999599,649768)$ & $649\,506\,443\,434$ & $3 \cdot 59 \cdot 3671$ & $(59,44)$ \\ \hline
\end{tabular}
\end{center}
\medskip

Note that there are two irregular pairs $(6449,4884)$ and $(6449,5830)$. But the first of them
disproves the suggested conjecture with minimal $l=4884$. The smallest index for which such an
exception occurs is $7\,610\,864$. This index is the smallest
of our calculated exceptions. For irregular pairs $(p,l)$ with
$p > 1\,000\,000$ we obtain index $m = (l-1) p+1 > 37\cdot 10^6$ for a possible exception,
because $37$ is the first irregular prime.

\section{Extending results to prime powers}
\setcounter{equation}{0}

In order to extend the results to irregular prime powers,
we need some further definitions and generalization. First,
the Kummer congruences generally state
for $r \geq 1$, $n, n'$ even, $p$ prime, and $p-1 \notdiv n$
\begin{equation} \label{eqn-kummer-congr-r}
   (1-p^{n-1}) \frac{B_n}{n} \equiv (1-p^{n'-1}) \frac{B_{n'}}{n'} \pmod{p^r}
\end{equation}
with $n \equiv n' \ \mod{\eulerphi(p^r)}$.
\medskip

The definition of irregular pairs can be extended to irregular prime powers
which was first introduced by the author \cite[Section 2.5]{kellner02irrpairord},
see also \cite{kellner04irrpowbn} for details and new results.
Here we will recall necessary facts.

\begin{defin} \label{def-ordn-irrpair}
A pair $(p,l)$ is called an \textsl{irregular pair of order $n$} if
$p^n \pdiv B_l/l$ with $2 \leq l < \eulerphi(p^n)$ and even $l$.
Let
\[
   \IRR_n := \{ (p, l) \sep p^n \pdiv B_l/l, \,\,
     2 \leq l < \eulerphi(p^n), \,\,  2 \pdiv l \}
\]
be the set of irregular pairs of order $n$. For a prime $p$
the \textsl{index} of irregular pairs of order $n$ is defined by
\[
    i_n(p) := \# \{ (p,l) \sep (p,l) \in \IRR_n \} \,.
\]
Let $(p,l) \in \IRR_n$ be an irregular pair of order $n$.
Let
\[
   (p,s_1,s_2, \ldots, s_n) \in \IRP_n \,, \quad
      l = \sum_{\nu=1}^{n} s_\nu \, \eulerphi( p^{\nu-1} )
\]
be the $p$-adic notation of $(p,l)$
with $0 \leq s_\nu < p$ for $\nu = 1,\ldots,n$ and $2 \pdiv s_1$, $2 \leq s_1 \leq p-3$.
The corresponding set will be denoted as $\IRP_n$. The pairs
$(p,l)$ and $(p,s_1,s_2, \ldots, s_n)$ will be called associated.
Define for an irregular pair $(p,l)$
\[
    \Delta_{(p,l)} \, \equiv \, p^{-1} \left(
       \frac{B_{l + \eulerphi(p)}}{l + \eulerphi(p)} - \frac{B_l}{l} \right) \pmod{p}
\]
with $0 \leq \Delta_{(p,l)} < p$.
\end{defin}
\bigskip

Note that this definition includes for $n=1$ the usual definition of irregular pairs
with $i(p) = i_1(p)$. By Kummer congruences \refeqn{eqn-kummer-congr-r} the interval
$[2,\eulerphi(p^n)-2]$ is given for irregular pairs of order $n$ if they exist.
Moreover, we have the property that if $(p,l) \in \IRR_n$ then
\begin{equation} \label{eqn-succ-pdiv-n}
   p^n \pdiv B_{l+k \eulerphi(p^n)}/(l+k \eulerphi(p^n))
\end{equation}
for all $k \in \NN_0$. Note that $(p,s_1,s_2, \ldots, s_n)$ is also called a \textsl{pair}
keeping in mind that $(s_1,s_2, \ldots, s_n)$ is the second parameter in a $p$-adic
manner. The main result of irregular pairs of higher order can be stated as follows,
see \cite[Theorem 3.1, p.~8]{kellner04irrpowbn}.

\begin{theorem} \label{theor-irrpair-main}
Let $(p,l_1)$ be an irregular pair. If $\Delta_{(p,l_1)} \neq 0$ then
for each $n > 1$ there exists exactly one irregular pair of order $n$
corresponding to $(p,l_1)$. Therefore a unique sequence
$(l_n)_{n \geq 1}$ resp. $(s_n)_{n \geq 1}$ exists with
\[
   (p,l_n) \in \IRR_n \quad \mbox{resp.} \quad (p,s_1,\ldots,s_n) \in \IRP_n \,.
\]
If $\Delta_{(p,l_{1,\nu)}} \neq 0$ for all $i(p)$ irregular pairs $(p,l_{1,\nu}) \in \IRR_1$, then
\[
   i(p) = i_2(p) = i_3(p) = \ldots \,\,.
\]
\end{theorem}

So far, no irregular pair $(p,l)$ with $\Delta_{(p,l)} = 0$ has been found
for $p < 12\,000\,000$ by calculations in \cite{irrprime12M}.
Because the case $\Delta_{(p,l)} = 0$ would
imply a strange behavior, it is conjectured that this will never happen.

\begin{theorem} \label{theor-Apr-small}
Let $r \geq 1$ be an integer.
Let $(p,l)$ be an irregular pair with $\Delta_{(p,l)} \neq 0$.
Then let $(p,l_r) \in \IRR_r$ resp. $(p,s_1,\ldots,s_r) \in \IRP_r$ be
the corresponding irregular pair of order $r$.
Define
\[
   A(p^r) = \min_k
     \left\{ m=l_r+k\eulerphi(p^r) \,\,\Big|\,\,
        \num \left( \frac{B_m}{m} \right) \Big/ \
           \num \left( \frac{B_m}{m(m-1)} \right) \ = \ p^r \right\} .
\]

Then $A(p^r)$ has only a solution if $(p,s_1,s_2,\ldots,s_r) = (p,l,l-1,\ldots,l-1)$ and
$l_r-1 = (l-1) p^{r-1} $. Furthermore $A(p^r) = (l_r-1)p+1 = (l-1)p^r+1$ is valid and
has smallest possible value if and only if one of the following cases holds
\begin{enumerate}
\item $l-1$ has no irregular prime factors.
\item If $q$ is an irregular prime divisor of $l-1$, then
all irregular pairs $(q,l')$ must satisfy
\[
   (l-1)p^r \ \not\equiv \ l'-1 \pmod{q-1} \,.
\]
\end{enumerate}
\end{theorem}

\begin{lemma} \label{lem-p-adic}
Let $n \geq 1$ and $s_1,\ldots,s_{n+1}$ be integers with $0 \leq s_\nu < p$ for all
$\nu=1,\ldots,n+1$. If
\[
   \sum_{\nu=1}^n s_\nu \, \eulerphi(p^{\nu-1}) = s_{n+1} \, p^{n-1} \,,
\]
then $s_1=s_2=\ldots=s_{n+1}$.
\end{lemma}

\begin{proof} Reordering terms yields
\[
   0 = \sum_{\nu=1}^n s_\nu \, \eulerphi(p^{\nu-1}) - s_{n+1} \, p^{n-1}
     = \sum_{\nu=1}^n (s_\nu - s_{\nu+1}) \, p^{\nu-1}
\]
which deduces the result $p$-adically by induction.
\end{proof}

\begin{proofof}{Theorem \ref{theor-Apr-small}}
Case $r = 1$ is handled by Theorem \ref{theor-Ap-small},
because $(p,l)=(p,l_1)=(p,s_1)$.
For now let $r \geq 2$. First we will show the proposed formula for $A(p^r)$.
To solve
\[
    \num \left( \frac{B_m}{m} \right) \Big/ \
        \num \left( \frac{B_m}{m(m-1)} \right) \ = \ p^r \,,
\]
factor $m-1$ must have the form $m-1 = p^r c$ with some integer $c$. Then $m-1$ must reduce
the $p$-power of the second numerator in order that
\[
    \ord_p \, \num \left( \frac{B_m}{m} \right) = u \quad \mbox{and} \quad
      \ord_p \, \num \left( \frac{B_m}{m(m-1)} \right) = u-r
\]
is valid with some integer $u \geq r$ which is granted by irregular pair $(p,l_r)$ of order $r$.
Let $m'$ be the smallest possible value.
By \refeqn{eqn-succ-pdiv-n}
we have $m' = l_r+k\eulerphi(p^r)$ and $m'-1 = p^r c$ which yields
\begin{equation} \label{eqn-loc-Apr-1}
  l_r-1 + k p^{r-1}(p-1) = p^r c  \quad \mbox{and} \quad  l_r-1 \equiv 0 \pmod{p^{r-1}} .
\end{equation}
Keeping in mind that
$l_r = \sum_{\nu=1}^{r} s_\nu \eulerphi( p^{\nu-1} ) < \eulerphi( p^r )$,
we obtain
\[
   0 < l_r-1 = p^{r-1} \, t < p^{r-1} (p-1)
\]
with $0 < t < p-1$. Rewriting \refeqn{eqn-loc-Apr-1} we get
\[
  p^{r-1} \, t + k p^{r-1}(p-1) = p^r c  \quad \mbox{and} \quad
    k p^{r-1} \equiv t p^{r-1} \pmod{p^{r}}
\]
which provides $k \equiv t \ \mod{p}$ and finally $k=t$ as smallest value.
Note that $l=s_1$ and $2 \leq l \leq p-3$. Now, using Lemma \ref{lem-p-adic}
with $l_r-1 = t p^{r-1}$ yields $s_1-1=s2=\ldots=s_r=t$.
Thus, we derive the following conditions
\[
   (p,s_1,s_2,\ldots,s_r) = (p,l,l-1,\ldots,l-1) \quad \mbox{and} \quad
     l_r-1 = (l-1) p^{r-1} \,.
\]
After all, we obtain
\[
   A(p^r)-1 = m'-1 = l_r - 1 + (l-1) \eulerphi( p^r ) = (l-1) p^r = (l_r-1) p \,.
\]
To avoid that an irregular prime divisor $q$ of the remaining factor $l-1$ of $m'-1$
divides $B_{m'}/m'$, we must have
\[
   m' \not\equiv \ l' \pmod{q-1}
\]
for all irregular pairs $(q,l')$. Then $A(p^r)$ is valid with the derived value.
\end{proofof}

\begin{corl} \label{corl-Apr-u}
Let $(p,l)$ be an irregular pair with $\Delta_{(p,l)} \neq 0$.
Let $r \geq 2$ be an integer, $(p,s_1,\ldots,s_r) \in \IRP_r$,
and $A(p^r)$ be defined as in Theorem \ref{theor-Apr-small}.
Assume $(p,s_1,s_2,\ldots,s_r)$ $\neq$ $(p,l,l-1,\ldots,l-1)$ then
$A(p^u)$ related to $(p,l)$ has no solution for all $u \ge r$.
\end{corl}

\begin{proof}
As a result of Theorem \ref{theor-irrpair-main}, if $\Delta_{(p,l)} \neq 0$
then a unique sequence $(s_\nu)_{\nu \ge 1}$ exists that
describes all irregular pairs of higher order related to $(p,l)$.
Then one has $(p,s_1,\ldots,s_r,\ldots,s_u)$ $\neq (p,l,l-1,\ldots,l-1)$
for all $u > r$.
\end{proof}
\medskip

The condition $(p,s_1,s_2,\ldots,s_r) = (p,l,l-1,\ldots,l-1)$
is a very strange condition. No such irregular pair $(p,s_1,s_2) \in \IRP_2$ of order two
with $s_2 = s_1-1$ has been found yet. For irregular primes $p<1000$ the
smallest difference $|s_1-s_2|$ is 4 which happens for the following elements
\[
   (353,186,190), \, (647,554,558) \in \IRP_2 \,.
\]
Therefore $A(p^r)$ has no solution for $p < 1000$ and $r \geq 2$.
Calculated irregular pairs of order 10 for $p<1000$ can be found in
\cite[Table A.3]{kellner04irrpowbn}.

\begin{remark} \label{rem-delta-0}
Although the more complicated case $\Delta_{(p,l)} = 0$ should not happen, Theorem
\ref{theor-Apr-small} is also valid in that case. We only need
an irregular pair $(p,l_r) \in \IRR_r$ and its associated pair
$(p,s_1,\ldots,s_r) \in \IRP_r$ which are related to $(p,l)$.
Corollary \ref{corl-Apr-u} remains to be valid in a similar way.
A strong condition must hold that further irregular pairs of order $r+1$
related to $(p,l_r)$ exist.
In case of existence they all have the form
$(p,s_1,\ldots,s_r,t) \in \IRP_{r+1}$ with $0 \leq t < p$,
see \cite[Theorem 3.2, p.~8]{kellner04irrpowbn}.
\end{remark}

\section{The composite case}
\label{sect-comp-case}
\setcounter{equation}{0}

For completeness we will examine the composite case. For now, we will recognize composite
integers $c$
\[
    c = \prod_{\nu=1}^n p_\nu^{e_\nu}
\]
having only irregular primes $p_\nu$ in its factorization with $n > 1$.
Therefore, $p$ will only denote irregular primes.
To determine the minimal index of the composite case, define
\[
   \Lambda(c) = \min_m
     \left\{ m \,\,\Big|\,\,
        \num \left( \frac{B_m}{m} \right) \Big/ \
           \num \left( \frac{B_m}{m(m-1)} \right) \equiv 0 \pmod{c} \right\} ,
\]
in case of no solution define $\Lambda(c) = \infty$.
Then, by Theorem \ref{theor-Ap-small}, we always have
\[
   \Lambda(p) = \min_{(p,l) \in \IRR_1} \ (l-1)p+1 \,.
\]
Theorem \ref{theor-Apr-small} asserts for $r \geq 2$
\[
   \Lambda(p^r) = \min_{(p,l,l-1,\ldots,l-1) \in \IRP_r} \ (l-1)p^r+1 \,,
\]
but there is no solution for $p < 1000$. Note that $m=12$ is the smallest index for
which $\num(B_m/m) > 1$. Hence, for $p > 1000$, $r \geq 2$, and $\Lambda(p^r) < \infty$,
we have a weak estimate
\begin{equation} \label{eqn-LamdaC-pr}
   \Lambda(p^r) > 11 \cdot 10^6 \,.
\end{equation}

\begin{lemma} \label{lem-lambda-cp}
Let $c = \prod_{\nu} p_\nu^{e_\nu}$ with irregular primes $p_\nu$. Then
\[
   \Lambda( c ) \, \geq \,
     \max_\nu \Lambda(p_\nu^{e_\nu}) \,.
\]
\end{lemma}

\begin{proof} Assume $\Lambda( c ) < \Lambda( p_\nu^{e_\nu} )$ for a fixed $\nu$.
But this contradicts the definition of $\Lambda$, because $p_\nu^{e_\nu} \pdiv c$.
The case of no solution is handled similarly.
\end{proof}
\medskip

Let $\mathcal{M}$ be the smallest index for which a composite number appears.
By our formerly calculated exceptions, we have an upper bound
\begin{equation} \label{eqn-min-lambda}
    \mathcal{M} \,=\, \min_c \Lambda(c) \,\leq\, 7\,610\,864 \,.
\end{equation}

Regarding estimate \refeqn{eqn-LamdaC-pr} for prime powers
above and using Lemma \ref{lem-lambda-cp},
for now, we only have to examine composite numbers which are squarefree.
Therefore, define the minimal value of $\Lambda$ for composite squarefree numbers
having $n \geq 2$ irregular prime factors by
\[
    \mathcal{M}_n = \min_{c=p_1 \cdots p_n} \Lambda(c) \,.
\]
Then, by definition we obviously have
\[
   \mathcal{M} = \mathcal{M}_2 \leq \mathcal{M}_3 \leq \ldots \,.
\]

For further results we need the well-known Chinese remainder theorem (CRT), s.
\cite[p. 34]{ireland90classical}, and its generalization.

\begin{theorem}[CRT]
Let $w_1,\ldots,w_n$ be positive integers which are pairwise relatively prime.
Define $W = \prod_{\nu=1}^n w_\nu$.
For a given system of simultaneous congruences
\[
   x \equiv a_\nu \pmod{w_\nu} \,, \qquad \nu=1,\ldots,n \,,
\]
there always exists a unique integer $x$ $\mod{W}$ with
\[
    x \equiv \sum_{\nu=1}^n a_\nu \, b_\nu \frac{W}{w_\nu} \pmod{W}
\]
and $b_\nu$ defined by
\[
    b_\nu \frac{W}{w_\nu} \equiv 1 \pmod{w_\nu} \,, \qquad \nu=1,\ldots,n \,.
\]
\end{theorem}

\begin{theorem}[CRT']
Let $w_1,\ldots,w_n$ be positive integers.
A system of simultaneous congruences
\[
   x \equiv a_\nu \pmod{w_\nu} \,, \qquad \nu=1,\ldots,n
\]
has a solution if and only if
\[
   a_i \equiv a_j \pmod{\gcd(w_i,w_j)}
\]
holds for all $i \neq j$.
Define $W = \lcm( w_1,\ldots,w_n )$, then $x$ has a unique solution $\mod{W}$.
\end{theorem}

To state our next theorem, we will introduce a new definition to characterize
a set of irregular pairs.

\begin{defin}
Irregular pairs $(p_1,l_1), \ldots, (p_n,l_n)$ are called \textsl{friendly} if
\[
    l_i \equiv l_j \pmod{\gcd(p_i-1,p_j-1)}
\]
is valid for all $i \neq j$. They are called \textsl{strong friendly} if, in addition,
\[
    p_i \not\equiv 1 \pmod{p_j} \qquad \mbox{or} \qquad
      (p_i,l_i) \equiv (1,1) \pmod{p_j}
\]
holds for all $i \neq j$.
\end{defin}
\medskip

For example, the irregular pairs (37,32), (59,44), (101,68) are strong friendly.
 \{(101,68), (607,592)\} and \{(131,22), (263,100)\} are sets of
friendly irregular pairs, but they are not strong friendly.

\begin{theorem} \label{theor-LamdaC}
Let $n \geq 2$ and $c=p_1 \cdots p_n$ be a composite number of
distinct irregular primes. Then $\Lambda(c)$ has only a solution if
there exists a set of strong friendly irregular pairs
$S=\{ (p_1,l_1), \ldots, (p_n,l_n) \}$. In case of existence there is
a unique integer $m_S$ with
\[
    c \, \le \, m_S-1 \, \le \, \lcm(c,p_1-1,\ldots,p_n-1)
\]
which simultaneously solves the congruences
\[
    m_S-1 \equiv p_\nu(l_\nu-1) \pmod{p_\nu(p_\nu-1)}\,,
      \quad \nu=1,\ldots,n \,.
\]
$\Lambda(c)$ is then given by
\[
    \Lambda(c) = \min_S m_S \,,
\]
whereas $S$ passes all such sets of strong friendly irregular pairs.
\end{theorem}

\begin{proof}
To derive conditions let $m$ be an integer solving
\[
      \num \left( \frac{B_m}{m} \right) \Big/ \
         \num \left( \frac{B_m}{m(m-1)} \right) \equiv 0 \pmod{c} \,.
\]
Thus, $c \pdiv B_m/m$ and $c \pdiv m-1$ provide
the existence of irregular pairs $(p_\nu,l_\nu)$ with
\begin{equation} \label{eqn-loc-LamdaC-1}
  \begin{array}{ll}
    m - 1 \,\equiv\, 0       & \mod{p_\nu} \\
    m - 1 \,\equiv\, l_\nu-1 & \mod{p_\nu-1} \\
  \end{array}
\end{equation}
for $\nu=1,\ldots,n$.
The system \refeqn{eqn-loc-LamdaC-1} of simultaneous congruences
has only a solution if conditions of CRT' are satisfied.
Therefore we have to recognize two cases
\begin{equation} \label{eqn-loc-LamdaC-2}
  \begin{array}{ll}
    l_i-1 \,\equiv\, l_j-1 & \mod{\gcd(p_i-1,p_j-1)} \\
    l_i-1 \,\equiv\, 0     & \mod{\gcd(p_i-1,p_j)} \\
  \end{array}
\end{equation}
which must be valid for all $i \neq j$.
The first congruence of \refeqn{eqn-loc-LamdaC-2} implies
that all considered irregular pairs must be friendly.
Additionally by the second congruence they must be strong friendly.
This property must hold for a solution and defines set $S$.
Combining \refeqn{eqn-loc-LamdaC-1} by CRT, we get
\begin{equation} \label{eqn-loc-LamdaC-3}
    m-1 \equiv p_\nu(l_\nu-1) \pmod{p_\nu(p_\nu-1)} \,, \quad
      \nu=1,\ldots,n \,.
\end{equation}
Let $W=\lcm(p_1(p_1-1),\ldots,p_n(p_n-1))$,
then system \refeqn{eqn-loc-LamdaC-1} resp.
\refeqn{eqn-loc-LamdaC-3} has a unique solution $\mod{W}$ by CRT'
and given set $S$.
Taking $1,\ldots,W$ as residue classes, we obtain a
minimal solution $m_S-1$ with the desired properties.
If $i(p_\nu) \geq 2$ holds for one index $\nu$, then probably other
sets $S$ can exist corresponding to irregular primes $p_1,\ldots,p_n$.
Therefore all such sets must be considered to get
\[
    \Lambda(c) = \min_S m_S \,. \vspace*{-4ex}
\]
\end{proof}
\medskip

Theorem \ref{theor-LamdaC} implies the following easy algorithm.

\begin{algo}
Let $n \geq 2$, $U$ be integers. Given an existing upper bound $U$ of $\mathcal{M}_n$,
define $u = \lfloor U^{1/n} \rfloor$. Otherwise set $U=u=\infty$.
Consider irregular primes
\begin{equation} \label{eqn-loc-algo-1}
  p_1 < \ldots < p_n \quad \mbox{with} \quad p_1 \cdots p_n < U \,, \quad p_1 < u \,.
\end{equation}
Start with smallest primes. For each tuple of primes do
\begin{itemize}
\item Step 1.
  Check for sets $S=\{ (p_1,l_1), \ldots, (p_n,l_n) \}$ of strong friendly irregular pairs.
  For each existing set $S$ calculate $m_S$ using Theorem \ref{theor-LamdaC}.
  Let $m=\min_S \, m_S$. If $m < U$ update $U \leftarrow m$ and $u$.
\item Step 2.
  If possible go to next primes satisfying \refeqn{eqn-loc-algo-1},
  otherwise stop with $\mathcal{M}_n=U$.
\end{itemize}
\end{algo}
\bigskip

Starting with $n=2$ and $U = 7\,610\,864$ yields $\mathcal{M}_2 = 107\,430$ with
$c=103 \cdot 149$.
Thus $\mathcal{M} = 107\,430$ is the smallest index for which a composite value
occurs. The result for $n=3$ is a quite large number with $\mathcal{M}_3 = 3\,754\,314\,782$,
see table below.
To check this result, irregular pairs $(p,l)$ up to $p < 2\,000\,000$ must be considered
for the first small primes.

\begin{center}
\begin{tabular}{|c|c|c|c|} \hline
$n$ & $S$ & $U$ & $u$ \\ \hline\hline
2 & $\{(37,32), \,(59,44)\}$    & $272\,876$ & 522 \\ \hline
2 & $\{(103,24), \,(149,130)\}$ & $107\,430$ & 327 \\ \hline\hline
3 & $\{(37,32), \,(59,44), \, (101,68)\}$ & 3\,979\,497\,668 & 1584 \\ \hline
3 & $\{(157,62), \,(401,382), \, (1217,1118)\}$ & 3\,754\,314\,782 & 1554 \\ \hline
\end{tabular}
\end{center}
\medskip

All results were calculated by several C++ programs and finally checked with Mathematica.

\section{A connection with Iwasawa theory}
\setcounter{equation}{0}

In Section \ref{sect-comp-case} we have seen that
Theorem \ref{theor-Apr-small} asserts for $r \geq 2$
\[
   \Lambda(p^r) = \min_{(p,l,l-1,\ldots,l-1) \in \IRP_r} \ (l-1)p^r+1 \,,
\]
noting that there is no solution for $p < 1000$. For a solution with $r \geq 2$
we basically need the existence of an irregular pair $(p,l,l-1) \in \IRP_2$ of order two.
\medskip

Now, the remarkable fact is that conditions $\Delta_{(p,l)} \neq 0$
and $(p,l,l-1) \notin \IRP_2$ play an important role in Iwasawa theory of
cyclotomic fields over $\QQ$, see \cite[Section 6]{kellner04irrpowbn}.
Here we give a brief summary.
\medskip

Let $\QQ(\mu_{p^n})$ be the cyclotomic field and  $\QQ(\mu_{p^n})^+$
its maximal real subfield with $\mu_{p^n}$ as the set of $p^n$-th roots of unity
Define the class number $h_p = h(\QQ(\mu_p))$ and its factoring $h_p = h_p^- \, h_p^+$
with $h_p^+ = h(\QQ(\mu_{p})^+)$ and $h_p^-$ as
the relative class number introduced by Kummer.
For details of the following theorem,
see \cite[Corollary 10.17, p.~202]{washington97cyclo}.
Note that conditions (2) and (3) are equivalently exchanged by our definitions.

\begin{theorem}
Let $p$ be an irregular prime. Assume the following conditions
\begin{enumerate}
\item The conjecture of Kummer--Vandiver holds: $p \notdiv h_p^+$
\item The $\Delta$-Conjecture holds: $\Delta_{(p,l)} \neq 0$
\item A special irregular pair of order two does not exist: $(p,l,l-1) \notin \IRP_2$
\end{enumerate}
for all irregular pairs $(p,l)$. Then
$\ord_p h( \QQ( \mu_{p^n} ) ) \, = \, i(p) \, n $ is valid for all $n \geq 1$.
\end{theorem}

Buhler, Crandall, Ernvall, Mets\"ankyl\"a, and Shokrollahi \cite{irrprime12M}
have calculated not only irregular pairs, but also associated
cyclotomic invariants up to $p < 12\,000\,000$. These calculations ensure that
no irregular pair $(p,l,l-1) \in \IRP_2$ exists in that range.
\medskip

Therefore we have a much stronger estimate than \refeqn{eqn-LamdaC-pr}
\[
   \Lambda(p^r) > 1.729 \cdot 10^{15}
\]
which can be obviously improved by choosing a greater value $l > 12$ examining the
numerators of the first divided Bernoulli numbers $B_m/m$.

\subsection*{Acknowledgement}
The author wishes to thank Tony D. Noe for advising the problem and Sequence \textbf{A092291}.
\bigskip

Bernd C. Kellner \\
{\small address: Reitstallstr. 7, 37073 G\"ottingen, Germany \\
email: bk@bernoulli.org}

\bibliographystyle{alpha}
\bibliography{conj_bib}

\end{document}